\documentclass[11pt,reqno]{amsart}
\usepackage[T1]{fontenc}
\usepackage{amsfonts}
\usepackage{amssymb}
\usepackage{mathrsfs}
\usepackage[latin1]{inputenc}
\usepackage{amsmath}
\usepackage{paralist}
\usepackage[english]{babel}

\newtheorem{theorem}{Theorem}[section]
\newtheorem{lemma}[theorem]{Lemma}

\newtheorem{remark}[theorem]{Remark}
\newtheorem{prop}[theorem]{Proposition}

\newtheorem{corollary}[theorem]{Corollary}

\numberwithin{equation}{section}
\newcommand{\R}{{\mathbb R}}

\newcommand{\N}{{\mathbb N}}

\newcommand{\Pb}{P_t^{\gamma,b}}
\newcommand{\Po}{P_t^{\gamma,0}}
\newcommand{\pb}{p^{\gamma, b}(x,y,t)}
\newcommand{\Ab}{A^{\gamma,b}}
\begin{document}

\title{One-dimensional degenerate diffusion  operators}

\author{Angela A. Albanese, Elisabetta M. Mangino}



\address{ Angela A. Albanese\\
Dipartimento di Matematica e Fisica ``E. De Giorgi''\\
Universit\`a del Salento\\
I-73100 Lecce, Italy}
\email{angela.albanese@unisalento.it}

\address{ Elisabetta M. Mangino\\
Dipartimento di Matematica e Fisica  ``E. De Giorgi''\\
Universit\`a del Salento\\
I-73100 Lecce, Italy} \email{elisabetta.mangino@unisalento.it}

\begin{abstract}
The aim of this paper is to  present some results about generation, sectoriality and gradient estimates both for the semigroup and for the resolvent  
 of   suitable realizations of the operators
\[ \Ab u(x)=\gamma xu''(x) + b u'(x), \]
with constants $\gamma >0$ and $b\geq 0$, in the space $C([0,\infty])$.
\end{abstract}

\maketitle
The motivation for this paper comes from  investigations on the  analiticity in the space of continuous functions on the $d$-dimensional canonical simplex $S^d$ of the semigroup generated by the  multi-dimensional Fleming-Viot operator (also known as Kimura operator or Wright-Fischer operator, see \cite{E,EK,EK1,FV,S1,S2})

\begin{equation}\label{e.operator}
Au(x)=\frac{1}{2}\gamma(x)\sum_{i,j=1}^dx_i(\delta_{ij}-x_j)\partial_{x_ix_j}^2u(x)+\sum_{i=1}^db_i(x)\partial_{x_i}u(x),
\end{equation}
where $b=(b_1, \dots,b_d)$ is a continuous inward pointing drift and $\gamma$ a strictly positive continuous function on $S^d$.
The operator \eqref{e.operator} arises in the theory of Fleming--Viot processes as the generator of  a Markov $C_0$--semigroup defined on $C(S_d)$. Fleming--Viot processes are measure--valued processes that can be viewed as diffusion approximations of empirical processes associated with some classes of discrete time Markov chains in population genetics. We refer to \cite{EK,EK1,FV} for more details on the  topic.

If  $b=0$, it has been proved in  \cite{AM} that the closure of  $(A,C^2(S^d))$ generates a bounded analytic semigroup, but  to extend the result to the case of a non-vanishing drift, it is needed  a careful estimate of the constants appearing in  the study of the sectoriality of the one-dimensional operator 
\begin{equation}\label{eq:fv}
Au(x)=\gamma(x)x(1-x)u''(x) + b(x)u'(x),\qquad x\in [0,1],\end{equation}
where $\gamma$ is a continuous strictly positive function and $b$  is a continuous function such that $b(0)\geq0$ and $b(1)\leq 0$.
As already pointed out by Feller in the fifties (see \cite{F1,F2}, see also \cite{T,CT,CMPR}), these degenerate operators generate positive and contractive semigroups in $C([0,1])$ if suitable boundary conditions are added. 

The aim of this paper is to  present some results about generation, sectoriality and gradient estimates for the resolvent  
 of   suitable realizations of the operators
\begin{equation}\label{eq:1} \Ab u(x)=\gamma xu''(x) + b u'(x), \end{equation}
with constants  $\gamma >0$ and $b\geq 0$, in the space $C([0,\infty])$, because they model in $0$
the behaviour of the operators (\ref{eq:fv})  near the end points $0$ and $1$.

To this end  we review some  results established mainly in \cite{BP, Met, EM, EM-1},  and we propose them in a unified way. 
We point out that the works  \cite{BP, EM, EM-1} are addressed mainly to the study of the operators (\ref{eq:fv}), (\ref{eq:1}) in H\"older continuous function spaces,
while we will mainly focus on spaces of continuous functions.
Moreover, several proofs are different from those of the cited papers  and could be of independent interest.

Precisely, we start considering the explicit expression of the kernel $\pb$, given, e.g., in \cite{EM,BP}, of the solution operator $\Pb$ for the equation $\partial_t-\Ab$.
After proving some estimates for $\pb$, we show that $(\Pb)_{t\geq 0}$ is a $C_0$-semigroup  in $C([0,\infty])$ and that its infinitesimal generator is $\Ab$ endowed with the domain
\begin{eqnarray*}
D(A^{\gamma,0})&=&\{ u\in C([0,\infty]) \cap C^2(]0,\infty[)\, \mid\,  \lim_{x\to 0^+} A^{\gamma,0} u(x) =0, \\
& & \lim_{x\to +\infty}A^{\gamma,0}u(x)=0\}, \qquad {\rm if}\ b=0,\\
 D(\Ab)&=&\{ u\in C^1([0,\infty[)\cap C^2(]0,\infty[)\cap C([0,\infty])\,\mid\, \\
& &\lim_{x\to 0^+} xu''(x)=0,\ \lim_{x\to + \infty} \Ab u(x)=0 \}, \qquad {\rm if}\ b>0.
\end{eqnarray*}
Moreover, we prove that the space of $C^2$-functions on $[0,\infty[$ that are constant in a neighbourhood of $\infty$ is  a core for  $(\Pb)_{t\geq 0}$.

At this point, the analiticity of $(\Pb)_{t\geq 0}$ in $C([0,\infty])$ follows immediately from  the results in \cite{Met, CM}, but with a careful analysis we also prove  that,  for any $B>0$ and $\gamma_0>0$ fixed, there exists a constant  $C=C(B,\gamma_0)>0$ such that,
for every $b\in [0,B]$ and  $\gamma\geq \gamma_0$,  
\[ 
||t\Ab\Pb || \leq C(B, \gamma_0),\quad t\geq 0, 
\] 
that is,  the analiticity constant is uniform in bounded intervals $[0,B]$ and in half-lines $[\gamma_0,\infty[$ with $\gamma_0>0$.

We also get pointwise gradient estimates both for the semigroup and for  the resolvent $R(\lambda, \Ab)$ and, in the case $b>0$, we prove that 
$\partial_x R(\lambda, \Ab)$ is a continuous operator from $C([0,\infty])$ into itself and give an estimate of the operator norm. 

The results presented in this paper play an important role in \cite{AM1} to show the analyticity in spaces of continuous functions of the semigroup generated by some degenerate diffusion operators defined on domains of $\R^d$ with corners  like \eqref{e.operator}. For further results on regularity in weighted $L^p$ spaces of the semigroup generated by some classes  of   operators  of type \eqref{e.operator} we refer to \cite{AM-0,ACM} and the references therein.

\bigskip

{\bf Notation.} We will denote by  $C_\mathbf{b}([0,\infty[)$ the space of continuous bounded functions on $[0,\infty[$ and by $C([0,\infty])$ the Banach space of continuous functions on $[0,\infty[$ converging  at infinity,  endowed with the sup-norm $||\cdot||_\infty$.
Analogously, for every $k\in\N$, $C^k([0,\infty])$ stands for the space of functions $u\in C([0,\infty])$ with derivatives up to order $k$ that have finite limits at $\infty$. Finally $C_c^k([0,\infty[)$ denotes the subspace of $C^k([0,\infty[)$ of functions with compact support and $C_0([0,\infty[)$ denotes the space of continuous functions on $[0,\infty[$ vanishing at $\infty$.

\section{Preliminary results}
 
\begin{lemma}\label{le:B} Let $A>0$.  Then there exists $C=C(A)>0$ such that, for every $0<a\leq A$ and $s>0$,
\begin{equation}
\sum_{m=0}^\infty \frac{s^m}{m!\Gamma(m+a)} \leq C \frac{e^{2\sqrt s}}{s^{\frac{a}{2}-\frac 1 4}}\left(1+\frac{e^\frac{C}{2\sqrt s}}{\sqrt{s}}\right).
\end{equation}
\end{lemma}

{\sc Proof.} Recall that the  Bessel modified function of the first kind and parameter $\nu\in\R$ is defined by the formula
\begin{equation} \label{Bessel} 
I_\nu(x)= \sum_{m=0}^\infty \frac{\left(\frac x 2\right)^{\nu+2m}}{m!\Gamma(m+\nu+1)}, \quad x>0,
\end{equation}
and so, for $\nu=a-1$ we have 
\[ 
\sum_{m=0}^\infty \frac{s^m}{m!\Gamma(m+a)}= \frac{1}{s^{\frac{a-1}{2}}} I_{a-1}(2\sqrt s), \quad s>0.
\]
Moreover, by \cite[(7.16)]{O}, for every $x>0$, the following estimate holds
\begin{eqnarray*}
|I_\nu(x)|&\leq& \frac {1}{(2\pi x)^{\frac 1 2}} \left[e^{ x}\left( |4\nu^2-1| + \pi^\frac 1 2 
e^\frac{\left|\nu^2-\frac 1 4\right|}{x} \frac{|(4\nu^2-1)(4\nu^2-9)|}{x}\right)  \right.\\
 &+& \left. e^{- x}\left( |4\nu^2-1| + 2e^\frac{\left|\nu^2-\frac 1 4\right|}{x} \frac{|4\nu^2-1|}{x} \right) \right] \leq 
\frac {C}{ x^{\frac 1 2}} e^{ x}\left(1+ e^{\frac{C}{x}}\frac 1 x\right),
\end{eqnarray*}
where $C=\max\{ \sqrt \pi |4\nu^2-1||4\nu^2-9|, 2| 4\nu^2-1|\}$. Then
 the assertion follows by applying the previous estimate with $\nu=a-1$, $C=\max_{[0,A]}\{\sqrt \pi |4a^2-8a+3| |4a^2-8a- 5|, 2|4a^2-8a +3|\}$ and $x=2\sqrt{s}$. \qed

\begin{lemma}\label{lem:2}
For every $\delta>0$ and for every $s>0$
\[ 
\sum_{m=0}^\infty |m-s| \frac{s^m}{(m+1)!} \leq \delta\frac{e^s-1}{s} + \frac 1 \delta \left(e^s-2-s + \frac{e^s-1}{s}\right).
\]
In particular, we have
\begin{eqnarray}
\sum_{m=0}^\infty |m-s| \frac{s^m}{(m+1)!} = O(\sqrt{s})\ \ as\ s\to 0^+\\
e^{-s} \sum_{m=0}^\infty |m-s| \frac{s^m}{(m+1)!} =O(\frac{1}{\sqrt{s}})\ \ \ as\ s\to\infty.
\end{eqnarray}
\end{lemma}

{\sc Proof.} Fix any $\delta>0$. Then, for every $s>0$, we have
\begin{eqnarray*}
& &\sum_{m=0}^\infty |m-s| \frac{s^m}{(m+1)!} = \sum_{|m-s|\leq \delta} |m-s| \frac{s^m}{(m+1)!} + \sum_{|m-s|> \delta} |m-s| \frac{s^m}{(m+1)!} \\
&\leq& \delta \sum_{m=0}^\infty \frac{s^m}{(m+1)!}  + \frac 1 \delta \sum_{m=0}^\infty (m-s)^2 \frac{s^m}{(m+1)!} \\
&=&  \delta \frac{e^s-1}{s} + \frac 1 \delta\left( \sum_{m=1}^\infty  \frac{m^2s^m}{(m+1)!} + s\sum_{m=0}^\infty \frac{s^{m+1}}{(m+1)!} - 2s \sum_{m=0}^\infty  \frac{m s^m}{(m+1)!}\right) \\
 &=& \delta \frac{e^s-1}{s} + \frac 1 \delta\left( \sum_{m=1}^\infty  \frac{m(m+1)s^m}{(m+1)!} + s(e^s-1)- 2s \sum_{m=0}^\infty  \frac{m s^m}{(m+1)!}- \sum_{m=1}^\infty \frac{m s^m}{(m+1)!}\right) \\
 &=& \delta\frac{e^s-1}{s} + \\
& & +\frac 1 \delta\left[ s\sum_{m=1}^\infty \frac{s^{m-1}}{(m-1)!} + s(e^s-1) -
 (2s+1)\left( \sum_{m=1}^\infty \frac{(m+1)s^m}{(m+1)!}
- \sum_{m=1}^\infty \frac{s^m}{(m+1)!}\right)\right]\\
 &=&\delta \frac{e^s-1}{s} + \frac 1 \delta\left[ se^s+ s(e^s-1) - (2s+1)\left( \sum_{m=1}^\infty \frac{s^m}{m!} - \frac 1 s \sum_{m=1}^\infty \frac{s^{m+1}}{(m+1)!}\right)\right]\\
 &=&\delta \frac{e^s-1}{s} + \frac 1 \delta \left(e^s-2-s+\frac{e^s-1}{s}\right).
 \end{eqnarray*}
 Choosing $\delta = \sqrt{s}$, we get that 
 \[ 
\sum_{m=0}^\infty |m-s| \frac{s^m}{(m+1)!} \leq \frac{e^s-1}{\sqrt{s}} + \frac{1}{\sqrt s}\left(e^s-2-s+\frac{e^s-1}{s}\right).\]
 Since 
 \[
 \frac {1}{\sqrt s}\left(e^s-2-s+\frac{e^s-1}{s}\right)\sim \frac 1 2 \sqrt{s}\quad {\rm and}\quad \frac{e^s-1}{\sqrt{s}}\sim \sqrt{s}\ {\rm as}\ s\to 0,
\]
 we easily deduce the assertion for $s\to 0^+$, while the behaviour as $s\to\infty$ follows observing that 
 \[ 
 \frac{e^s-1}{\sqrt{s}} + \frac{1}{ \sqrt{s}} \left(e^s-2-s+\frac{e^s-1}{s}\right) \sim 2\frac{e^s}{\sqrt s}\ {\rm as}\ s\to\infty. \qquad \qed\]

\begin{lemma}\label{le:gamma} For every $s>0$ 
\[ 
\frac{1}{\Gamma(s)}\int_0^\infty e^{-z}z^{s-1}|z-s| dz \leq 2\sqrt{s}.
\]
\end{lemma}

{\sc Proof.} For every $\delta>0$ and $s>0$
\begin{eqnarray*}
& &\frac{1}{\Gamma(s)} \int_0^\infty e^{-z}z^{s -1}\left|z-s\right|dz \\
&\leq& \frac{1}{\Gamma(s)} \int_{|z-s|<\delta} e^{-z}z^{s -1}\left|z-s\right|dz + \frac{1}{\Gamma(s)} \int_{|z-s|\geq \delta} e^{-z}z^{s -1}\left|z-s\right|dz\\
&\leq& \frac{1}{\Gamma(s)} \left(\delta \Gamma(s) + \frac{1}{\delta} \int_0^\infty e^{-z}z^{s -1}\left(z-s\right)^2dz\right)\\
&=&
 \frac{1}{\Gamma(s)} \left[\delta \Gamma(s) + \frac{1}{\delta}\left( \Gamma(s+2) -2s \Gamma (s+ 1) + s^2\Gamma(s)\right)\right]\\
&=&\frac{1}{\Gamma(s)} 
\left(\delta \Gamma(s) + \frac{1}{\delta} \Gamma(s+1)\right)=  \delta  + \frac s \delta.
\end{eqnarray*}
Then the thesis follows by choosing $\delta=\sqrt{s}$.\qed

\section{The semigroup $(P_t^{\gamma,b})_{t\geq 0}$}
Fix any $\gamma>0$ and $b\geq 0$. For every $x,y\geq 0$ and $t>0$ set
\begin{eqnarray*}
 & & p^{\gamma, b}(x,y,t)= (\gamma t)^{-\frac b \gamma}e^{-\frac{x+y}{\gamma t}}y^{\frac b \gamma -1} \sum_{m=0}^\infty \frac{1}{m!\Gamma (m+\frac b \gamma)}\frac{x^m y^m}{\gamma^{2m}t^{2m}} \qquad {\rm if}\ b>0\\
& & p^{\gamma,0 }(x,y,t)= e^{-\frac{x}{\gamma t}}\delta_0(y)+(\gamma t)^{-1}e^{-\frac{x+y}{\gamma t} }\sum_{m=0}^\infty \frac{1}{(m+1)!m!}\frac{x^{m+1} y^m}{\gamma^{2m+1}t^{2m+1}}\qquad {\rm if}\ b=0,
\end{eqnarray*}
where $\delta_0$ is a point mass at $0$.
Thanks to (\ref{Bessel}) the kernels can  be  written as 
\begin{eqnarray} & &p^{\gamma,b}(x,y,t)= \frac{1}{\gamma t} \left(\frac x y \right)^{\frac{1-b/\gamma}{2}} e^{-\frac{x+y}{\gamma t}} I_{b/\gamma -1}\left( \frac{2\sqrt{xy}}{\gamma t}\right)
\qquad {\rm if}\ b>0\\
& & p^{\gamma,0}(x,y,t)=e^{-\frac{x}{\gamma t}}\delta_0(y)+\frac{1}{\gamma t} \left(\frac x y \right)^{\frac{1}{2}} e^{-\frac{x+y}{\gamma t}} I_{1}\left( \frac{2\sqrt{xy}}{\gamma t}\right)\qquad {\rm if}\ b=0.
\end{eqnarray}
Therefore, we easily obtain that
\begin{equation}\label{e.ident-one}
 \int_0^\infty p^{\gamma, b}(x,y,t)dy=1,\qquad   x\geq 0,\ t>0.
\end{equation}
Moreover, in \cite[Corollary 9]{EM} (see also \cite{BP} for  a stochastic approach), it is shown that 
\[ 
p^{\gamma, b}(x,y,t+s)=\int_0^\infty p^{\gamma,b}(x,z,t) p^{\gamma,b}(z,y,s)dz\qquad \ s,\, t>0,\  x,\, y\geq 0.
\]
The  properties above ensure that, for each $t>0$, the operator $P_t^{\gamma,b}\colon C_\mathbf{b}([0,\infty[)\to C_\mathbf{b}([0,\infty[)$ defined by 
\begin{equation}\label{e.defsem}
P_t^{\gamma,b} f(x)= \int_0^\infty p^{\gamma,b} (x,y,t)f(y)dy, \qquad f\in C_\mathbf{b}([0,\infty[),\ x\geq 0,
\end{equation} 
is well-defined and continuous and that the family $(\Pb)_{t\geq 0}$ is a contraction semigroup in $C_\mathbf{b}([0,\infty[)$  (here, $P_0^{\gamma,b}:=I$).
We  also observe that, for every $f\in C_\mathbf{b}([0,\infty[)$ and $x\geq 0$, we have
\begin{equation}\label{e.rapsem} 
\Pb f(x)=e^{-\frac{x}{\gamma t}}\sum_{m=0}^\infty \left(\frac{x}{\gamma t}\right)^m \frac{1}{m!\Gamma (m+\frac b \gamma)} \int_0^\infty e^{-z}z^{m+\frac b \gamma -1} f(\gamma z t)dz \qquad {\rm if}\ b>0,
\end{equation}
while 
\begin{equation}\label{e.rapsem0}
 P^{\gamma,0}_tf(x)=e^{-\frac{x}{\gamma t}}f(0)+e^{-\frac{x}{\gamma t}}\sum_{m=0}^\infty \left(\frac{x}{\gamma t}\right)^{m+1} \frac{1}{m!(m+1)!} \int_0^\infty e^{-z}z^{m} f(\gamma z t)dz \ {\rm if}\ b=0.
\end{equation}
Since $\int_0^\infty \pb y^k dy$ is  convergent for every $k\in\N$,   with an abuse of notation we can set $\Pb(f)=\int_0^\infty \pb f(y)dy$ for every continuous function $f$ with polynomial growth at infinity. 

In the following,   for each $x,\, y>0$ we set $\tau_x(y):=(y-x)$.

\begin{lemma}\label{le.semb} Let $b>0$ and $\gamma>0$. Then,  for every $x,t>0$ and $k\in\N$, the following properties hold:
\begin{itemize}
\item[\rm (1)] $\Pb(\tau_x)(x)=b t$,
\item[\rm (2)] $\Pb(\tau_x^{k+1})(x)=-x\Pb(\tau_x^k)(x)+ btP_t^{\gamma,b+\gamma}(\tau_x^k)(x) +  xP_t^{\gamma,b+2\gamma}(\tau_x^k)(x)$.
\end{itemize}
In particular, for every $x,t>0$, we have
\begin{eqnarray}
\label{eq:2} & &\Pb(\tau_x^2)(x)=2\gamma tx + t^2b(b+\gamma),\\
\label{eq:3n} & & \Pb(|\tau_x|)(x)  \leq \sqrt{ 2\gamma t x + t^2 b (b+\gamma)}.
\end{eqnarray}
\end{lemma}

{\sc Proof.} Let $b>0$. By \eqref{e.ident-one} and \eqref{e.rapsem} we obtain, for every $x,t>0$, that
\begin{eqnarray*}
\Pb(\tau_x)(x)&=& e^{-\frac{x}{\gamma t}}\sum_{m=0}^\infty \left(\frac{x}{\gamma t}\right)^m \frac{1}{m!\Gamma (m+\frac b \gamma)}\gamma t  \int_0^\infty e^{-z}z^{m+\frac b \gamma} dz - x\\
&=& e^{-\frac{x}{\gamma t}}\sum_{m=0}^\infty \left(\frac{x}{\gamma t}\right)^m \frac{m+b/\gamma}{m!}\gamma t-x \\
&=& x + b t -x = bt;
\end{eqnarray*}
hence, (1) is satisfied. 

Fixed any $k\in\N$, we have, for every $x,t>0$, that 
\begin{eqnarray*}
 & & \Pb(\tau_x^{k+1})(x)+x\Pb(\tau_x^k)(x)= \int_0^\infty \pb y(y-x)^k dy  \\
 &=& \gamma te^{-\frac{x}{\gamma t}}\sum_{m=0}^\infty \left(\frac{x}{\gamma t}\right)^m \frac{1}{m!\Gamma (m+\frac b \gamma)}  \int_0^\infty e^{-z}z^{m+\frac b \gamma}(\gamma t z- x)^k dz \\
 &=& \gamma te^{-\frac{x}{\gamma t}}\sum_{m=0}^\infty \left(\frac{x}{\gamma t}\right)^m \frac{m + \frac b \gamma}{m!\Gamma (m+\frac b \gamma+1)}  \int_0^\infty e^{-z}z^{m+\frac b \gamma}(\gamma t z- x)^k dz \\
&=& \gamma t e^{-\frac{x}{\gamma t}}\sum_{m=1}^\infty \left(\frac{x}{\gamma t}\right)^m \frac{1}{(m-1)!\Gamma (m+\frac b \gamma+1)}  \int_0^\infty e^{-z}z^{m+\frac b \gamma} (\gamma t z- x)^k dz \\ 
& &+\gamma t \frac b \gamma e^{-\frac{x}{\gamma t}}\sum_{m=0}^\infty \left(\frac{x}{\gamma t}\right)^m \frac{1}{m!\Gamma (m+\frac b \gamma + 1)}  \int_0^\infty e^{-z}z^{m+\frac b \gamma}(\gamma t z- x)^k dz\\
&=& x P_t^{\gamma,b+2\gamma}(\tau_x^k)(x) + btP_t^{ \gamma,b+\gamma}(\tau_x^k)(x)
\end{eqnarray*}
and so (2) is also satisfied.
As a consequence we deduce,  for every $x,t>0$, that 
\begin{eqnarray*}
\Pb(\tau_x^2)(x)&=&  -x \Pb(\tau_x)(x) + x P_t^{ \gamma,b+2\gamma,}(\tau_x)(x) + btP_t^{\gamma,b+\gamma}(\tau_x)(x)\\
&=& -xbt + x (b+2\gamma)t + bt(b+\gamma)t=2\gamma t x + t^2 b (b+\gamma).
\end{eqnarray*}
Finally, by applying H\"older inequality together with \eqref{e.ident-one} and  \eqref{eq:2}, we get, for every $x>0$, that  
\[ \int_0^\infty p^{\gamma,b}(x,y,t)|x-y|dy \leq 
 \left(\int_0^\infty  p^{\gamma,b}(x,y,t)|x-y|^2dy\right)^{\frac 1 2}= 
\sqrt{2\gamma t x + t^2 b (b+\gamma)}. \qquad\qed
\]

\begin{lemma}\label{le.sem0}
Let  $\gamma>0$. Then, for every $x,t>0$ the following properties hold:
\begin{itemize}
\item[\rm (1)] $P_t^{\gamma,0}(\tau_x)(x)=0$,
\item[\rm (2)]  $P_t^{\gamma,0}(\tau_x^2)(x)=2\gamma t x$.
\end{itemize}
In particular, for every $x,t>0$, we have
\begin{equation}\label{eq:sem0}
P_t^{\gamma,0}(|\tau_x|)(x)\leq \sqrt{2\gamma t x}.
\end{equation}
\end{lemma}

{\sc Proof.} By  \eqref{e.rapsem0} we obtain, for every $x,t>0$, that
\begin{eqnarray*}
P_t^{\gamma,0}(\tau_x)(x)&=&-xe^{-\frac{x}{\gamma t}}+e^{-\frac{x}{\gamma t}}\sum_{m=0}^\infty\left(\frac{x}{\gamma t}\right)^{m+1}\frac{1}{m!(m+1)!}\int_0^\infty e^{-z}z^m(\gamma t z -x)dz\\
&=& -xe^{-\frac{x}{\gamma t}}+\gamma t e^{-\frac{x}{\gamma t}}\sum_{m=0}^\infty\left(\frac{x}{\gamma t}\right)^{m+1}\frac{1}{m!}\\
& &-xe^{-\frac{x}{\gamma t}}\sum_{m=0}^\infty\left(\frac{x}{\gamma t}\right)^{m+1}\frac{1}{(m+1)!}\\
&=& -xe^{-\frac{x}{\gamma t}}+x-xe^{-\frac{x}{\gamma t}}(e^{\frac{x}{\gamma t}}-1)\\
&=& -xe^{-\frac{x}{\gamma t}}+x-x+xe^{-\frac{x}{\gamma t}}=0;
\end{eqnarray*}
hence, (1) is satisfied. Also, for every $x,t>0$, we have
\begin{eqnarray*}
P_t^{\gamma,0}(\tau_x^2)(x)&=& x^2e^{-\frac{x}{\gamma t}}+e^{-\frac{x}{\gamma t}}\sum_{m=0}^\infty\left(\frac{x}{\gamma t}\right)^{m+1}\frac{1}{m!(m+1)!}\int_0^\infty e^{-z}z^m(\gamma t z -x)^2dz\\
&=& x^2e^{-\frac{x}{\gamma t}}+\gamma^2 t^2e^{-\frac{x}{\gamma t}}\sum_{m=0}^\infty\left(\frac{x}{\gamma t}\right)^{m+1}\frac{(m+2)!}{m!(m+1)!}\\
& & -2\gamma t xe^{-\frac{x}{\gamma t}}\sum_{m=0}^\infty\left(\frac{x}{\gamma t}\right)^{m+1}\frac{(m+1)!}{m!(m+1)!}
 +x^2e^{-\frac{x}{\gamma t}}\sum_{m=0}^\infty\left(\frac{x}{\gamma t}\right)^{m+1}\frac{m!}{m!(m+1)!}\\
&=& x^2e^{-\frac{x}{\gamma t}}+ \gamma^2 t^2e^{-\frac{x}{\gamma t}}\frac{x^2}{\gamma^2 t^2}\sum_{m=0}^\infty\left(\frac{x}{\gamma t}\right)^{m}\frac{1}{m!}+2\gamma^2 t^2e^{-\frac{x}{\gamma t}}\frac{x}{\gamma t} \sum_{m=0}^\infty\left(\frac{x}{\gamma t}\right)^{m}\frac{1}{m!}\\
& & -2\gamma t x e^{-\frac{x}{\gamma t}}\frac{x}{\gamma t}e^{\frac{x}{\gamma t}}+x^2e^{-\frac{x}{\gamma t}}(e^{\frac{x}{\gamma t}}-1)\\
&=& x^2e^{-\frac{x}{\gamma t}}+x^2+2\gamma t x-2x^2+x^2-x^2e^{-\frac{x}{\gamma t}}=2\gamma tx
\end{eqnarray*}
and so (2) is satisfied. Finally, by applying H\"older inequality together with  \eqref{e.ident-one} and  the property (2) above one easily shows \eqref{eq:sem0}.\qed

\begin{lemma}\label{lim} Let $b\geq 0$ and $\gamma>0$. Then the following properties hold.
\begin{itemize} 
\item[\rm (1)] If $f\in C_\mathbf{b}([0,\infty[)$, then  
\[ 
\lim_{t\to 0^+}\Pb f=f 
\]
uniformly on compact subsets of  $[0,\infty[$.
\item[\rm (2)] If $b>0$,  then there exists a constant $C=C(b)>0$ such that, for every $f\in C_c(\R^+)$ with supp$(f)\subseteq [0,M]$ and $x, t>0$, we have
\begin{equation}\label{eq:lim1}
|\Pb f(x)|\leq C||f||_\infty e^{-\frac{x-2\sqrt{xM}}{\gamma t}} \left(\frac M x\right)^{\frac{b}{2\gamma} + \frac 1 4}  \sqrt{\frac{\gamma t}{M}} 
\left( 1+ e^{\frac{C\gamma t}{2\sqrt{xM}}} \frac{\gamma t}{\sqrt{xM}}\right). 
\end{equation}
If $b=0$, then 
 there exists a constant $C>0$ such that, for every $f\in C_c(\R^+)$ with supp$(f)\subseteq [0,M]$ and $x, t>0$, we have
\begin{equation}\label{eq:lim2}
|\P_t^{\gamma,0} f(x)|\leq C||f||_\infty e^{-\frac{x-2\sqrt{xM}}{\gamma t}} \left(\frac M x\right)^{\frac{1}{4}}  \sqrt{\frac{\gamma t}{M}} 
\left( 1+ e^{\frac{C\gamma t}{2\sqrt{xM}}} \frac{\gamma t}{\sqrt{xM}}\right). 
\end{equation}
Therefore, for every $t>0$, $\lim_{x\to\infty}\Pb f(x)=0$ and $\lim_{t\to 0^+}\Pb f(x)=f(x)$ uniformly on $[0, +\infty[$.
\end{itemize}
\end{lemma}

{\sc Proof.} 
(1) Let $f\in C_\mathbf{b}([0,\infty[)$ and let $M>0$. We prove that $\lim_{t\to 0}\Pb f(x)=f(x)$ uniformly in $[0,M]$. Indeed, let $\varepsilon >0$ and let $\delta>0$ be such that $|f(x)-f(y)|<\varepsilon$ whenever $x,y\in [0,M]$ satisfy $|x-y|<\delta$. Then, by  \eqref{eq:3n} and \eqref{eq:sem0}, we obtain, for every $x\in [0, M]$ and $t>0$, that
\begin{eqnarray*}
|\Pb f(x)-f(x)|&=&\left|\int_0^\infty p^{\gamma,b}(x,y,t) (f(y)-f(x))dy\right| \\
&\leq & \varepsilon \int_{|x-y|<\delta} p^{\gamma, b}(x,y,t)dy + 2||f||_\infty \int_{|x-y|\geq \delta} p^{\gamma,b}(x,y,t)dy\\
&\leq& \varepsilon + \frac{2||f||_\infty}{\delta}\int_0^\infty p^{\gamma,b}(x,y,t)|x-y| dy \\
&\leq &\varepsilon + 
\frac{2||f||_\infty}{\delta}\sqrt{2\gamma t x + t^2 b \left(  b+ \gamma\right)} \\
&\leq& \varepsilon + \frac{2||f||_\infty}{\delta}\sqrt{2\gamma t M + t^2 b \left(  b+ \gamma \right)}.
\end{eqnarray*}
 We now get immediately the assertion.

(2) Let $f\in C_c([0,\infty[)$ with supp$(f)\subseteq [0,M]$. If $b>0$, then, by \eqref{e.rapsem} and Lemma \ref{le:B}, we obtain,  for every $x,t>0$, that
\begin{eqnarray*}
|\Pb f(x)| &\leq & ||f||_\infty e^{-\frac{x}{\gamma t}}\sum_{m=0}^\infty \left( \frac{x}{\gamma t}\right)^m \frac{1}{m!\Gamma(m+b/\gamma)} \int_0^{\frac{M}{\gamma t}} e^{-z}z^{m+b/\gamma -1}dz\\
& \leq& ||f||_\infty e^{-\frac{x}{\gamma t}} \sum_{m=0}^\infty \left( \frac{x}{\gamma t}\right)^m \frac{1}{m!\Gamma(m+b/\gamma)} \frac{1}{m+b/\gamma}\left(\frac{M}{\gamma t}\right)^{m+b/\gamma}\\
&= & ||f||_\infty e^{-\frac{ x}{\gamma t}} \left(\frac{M}{\gamma t}\right)^{b/\gamma} \sum_{m=0}^\infty \left( \frac{xM}{\gamma^2 t^2}\right)^{m} \frac{1}{m!\Gamma(m+b/\gamma+1)}\\
&\leq & C||f||_\infty e^{-\frac{ x}{\gamma t}} \left(\frac{M}{\gamma t}\right)^{b/\gamma} e^{\frac{2\sqrt{xM}}{\gamma t}}\frac{1}{\left(\frac{xM}{\gamma^2 t^2}\right)^{\frac{b}{2\gamma}+\frac{1}{4}}}\left(1+e^{\frac{C\gamma t}{2\sqrt{xM}}}\frac{\gamma t}{\sqrt{xM}}\right)\\
&=& C||f||_\infty e^{-\frac{ x-2\sqrt{xM}}{\gamma t}} \sqrt{\frac{\gamma t}{M}}\left(\frac{M}{x}\right)^{\frac{b}{2\gamma}+\frac{1}{4}}\left(1+e^{\frac{C\gamma t}{2\sqrt{xM}}}\frac{\gamma t}{\sqrt{xM}}\right).
\end{eqnarray*}
If $b=0$, then, by \eqref{e.rapsem0} and Lemma \ref{le:B}, we obtain,  for every $x,t>0$, that
\begin{eqnarray*}
|P^{\gamma, 0}_tf(x)|&\leq & ||f||_\infty e^{-\frac{ x}{\gamma t}}+||f||_\infty e^{-\frac{ x}{\gamma t}}\sum_{m=0}^\infty\left(\frac{x}{\gamma t}\right)^{m+1}\frac{1}{m!(m+1)!}\int_0^{\frac{M}{\gamma t}} e^{-z}z^mdz\\
&\leq & ||f||_\infty e^{-\frac{ x}{\gamma t}}\left(1+\sum_{m=0}^\infty\left(\frac{x}{\gamma t}\right)^{m+1}\frac{1}{(m+1)!(m+1)!}\left(\frac{M}{\gamma t}\right)^{m+1}\right)\\
&=& ||f||_\infty e^{-\frac{ x}{\gamma t}}\sum_{m=0}^\infty \left(\frac{xM}{\gamma^2 t^2}\right)^m\frac{1}{m\Gamma(m+1)}\\
&\leq & C ||f||_\infty e^{-\frac{ x}{\gamma t}} e^{\frac{2\sqrt{xM}}{\gamma t}}\frac{1}{\left(\frac{xM}{\gamma^2 t^2}\right)^{\frac{1}{4}}}\left(1+e^{\frac{C\gamma t}{2\sqrt{xM}}}\frac{\gamma t}{\sqrt{xM}}\right)\\
&=& C ||f||_\infty e^{-\frac{ x-2\sqrt{xM}}{\gamma t}}\sqrt[4]{\frac M x}\sqrt{\frac{\gamma t}{M}}\left(1+e^{\frac{C\gamma t}{2\sqrt{xM}}}\frac{\gamma t}{\sqrt{xM}}\right).
\end{eqnarray*}
Now, if $f\in C_c([0,\infty[)$ with supp$(f)\subseteq [0,M]$, then \eqref{eq:lim1} (\eqref{eq:lim2} for $b=0$) clearly implies that $\lim_{x\to\infty}\Pb f(x)=0$ for every $t>0$. On the other hand, fixed any $\varepsilon>0$, by \eqref{eq:lim1} (\eqref{eq:lim2} for $b=0$) there exists $N>M$ so that $|\Pb f(x)|<\varepsilon/2$ for every $x>N$ and $0<t\leq 1$. By (1) there is also $\overline{t}\in ]0,1]$ for which $\max_{x\in [0,N]}|\Pb f(x)-f(x)|<\varepsilon/2$ for every $0<t<\overline{t}$. So, it follows, for every $0<t<\overline{t}$, that 
\[
\sup_{x\in [0,\infty[}|\Pb f(x)-f(x)| \leq \max_{x\in [0,N]}|\Pb f(x)-f(x)|+\sup_{x>N}|\Pb f(x)|<\varepsilon/2+\varepsilon/2=\varepsilon.
\]
This completes the proof.\qed

\begin{prop}\label{p:c0semig} Let $b\geq 0$ and $\gamma>0$. Then 
$(P_t^{\gamma,b})_{t\geq 0}$ is a $C_0$-semigroup in $C([0,\infty])$. 
\end{prop}

{\sc Proof.} 
Since each  operator $\Pb$ preserves constant functions and is contractive, it suffices to prove that $(P_t^{\gamma,b})_{t\geq 0}$ is a $C_0$-semigroup in $C_0([0,\infty[)$.  
So, we first observe that Lemma  \ref{lim}(2) and  the fact that $\Pb$ is a continuous linear  operator from $C_\mathbf{b}([0,\infty[)$ into itself ensure that
\[ 
\Pb(C_0([0,\infty[))= \Pb(\overline{C_c([0,\infty[)}) \subseteq \overline{\Pb(C_c([0,\infty[))} \subseteq C_0([0,\infty[).
\]
Hence,  $\Pb$ is a well-defined bounded linear operator from $C_0([0,\infty[)$ into itself. On the other hand, 
if  $f\in C_c([0,\infty[)$, then by Lemma \ref{lim}(2) we have $\lim_{t\to 0^+}\Pb f=f$ uniformly on $[0,+\infty[$. The density of $C_c([0,\infty[)$ in $C_0([0,\infty[$ and the contractivity of $(P_t^{\gamma,b})_{t\geq 0}$ imply that $\lim_{t\to 0^+}\Pb f=f$ uniformly on $[0,+\infty[$ for every $f\in C_0([0,\infty[)$. This completes the proof.\qed

\medskip

In the sequel, we denote by $(A, D(A))$ the generator of $(\Pb)_{t\geq 0}$ in $C([0,\infty])$ and set 
\begin{equation}\label{eq:core}
D:=\{C^2([0,\infty])\mid f \mbox{\ is constant in a neighbourhood of\ } +\infty\}.
\end{equation}

\begin{prop} Let $b\geq 0$, $\gamma>0$ and let $D$ be  defined according to \eqref{eq:core}. Then
$D\subseteq D(A)$ and  
\[Af(x)=\gamma x f''(x) + bf'(x)  \]
for every $f\in D$ and $x\geq 0$.
\end{prop}

{\sc Proof.} Fix any $f\in D$.  Since each $\Pb$ preserves constant functions, we can assume w.l.o.g. that  $f\in C^2([0,\infty[)$ with supp$(f)\subseteq [0,M]$. 

Now, let $b>0$. Then by \eqref{eq:lim1} we have, for every $x>9M$ and $t>0$, that
\begin{eqnarray*}
& & \left| \frac{\Pb f(x)-f(x)}{t}\right| \leq \\ 
&\leq &\frac{C||f||_\infty}{t} e^{-\frac{x-2\sqrt{xM}}{\gamma t}} \left(\frac M x\right)^{\frac{b}{2\gamma}+ \frac 1 4} \sqrt{\frac{\gamma t}{M}} \left( 1+ e^{\frac{C\gamma t}{2\sqrt{xM}}} \frac{\gamma t}{\sqrt{xM}}\right)\\
&\leq& \frac{C||f||_\infty}{t}e^{-\frac{3M}{\gamma t}}  \sqrt{\frac{\gamma t}{M}} \left( 1+ e^{\frac{C\gamma t}{M}} \frac{\gamma t}{M}\right).
\end{eqnarray*}
It follows that 
\begin{equation}\label{eq:gnt1}
\lim_{t\to 0^+}\frac{\Pb f(x)-f(x)}{t}=0= \gamma x f''(x)+bf'(x)\ {\rm uniformly\ on\ } [9M,\infty[.
\end{equation}
On the other hand, for every $x,y\geq 0$ we can write $f(y)=f(x)+f'(x)\tau_x(y) + \frac 1 2 f''(x)\tau_x^2(y)+ \omega(x,y)$, where ${\omega(x,y)}= \frac{f''(\xi)-f''(x)}{2}(x-y)^2$  with $\xi$ belonging to the interval having $x$ and $y$ as endpoints. 
Since $|\omega(x,y)|\leq ||f''||_\infty (x-y)^2$  for every $x, y\geq 0$, 
 by  (\ref{eq:2})  we obtain, for every $x, t>0$, that
\begin{eqnarray*}
& & \left| \frac{\Pb f(x)-f(x)}{t}- \gamma x f''(x)-bf'(x)\right| \\
&=& |\frac{1}{2}f''(x)tb(b+\gamma) +  \int_0^\infty \pb \omega(x,y) dy|  \\
&\leq & \frac{1}{2}t ||f''||_\infty b(b+\gamma) +||f''||_\infty \int_0^\infty\pb (y-x)^2dy \\
&= &  \frac{1}{2}t ||f''||_\infty b(b+\gamma) +  ||f''||_\infty \Pb(\tau_x^2)(x) \\
& = &\frac{1}{2}t ||f''||_\infty b(b+\gamma) +  ||f''||_\infty [2\gamma tx + t^2b(b+\gamma)].
\end{eqnarray*}
It follows that
\begin{equation}\label{eq:gnt2}
\lim_{t\to 0^+}\frac{\Pb f(x)-f(x)}{t}= \gamma x f''(x)+bf'(x)\ {\rm uniformly\ on\ compact\ subsets\ of\ } [0,\infty[. 
\end{equation}
So, by \eqref{eq:gnt1} and \eqref{eq:gnt2} we get that $\lim_{t\to 0^+}\frac{\Pb f(x)-f(x)}{t}=\gamma x f''(x)+bf'(x)$ in $C([0,\infty])$.

In case $b=0$ the proof is analogous via \eqref{eq:lim2} and Lemma \ref{le.sem0}(2).		\qed

\section{The infinitesimal generator $\Ab$}

We now consider the operator $\Ab u= \gamma xu''+bu'$,  with the domain $D(\Ab)$ defined in the introduction, i.e., 
\begin{eqnarray*}
D(A^{\gamma,0})&=&\{ u\in C([0,\infty]) \cap C^2(]0,\infty[)\, \mid\,  \lim_{x\to 0^+} A^{\gamma,0} u(x) =0, \\
   &  &\lim_{x\to +\infty}A^{\gamma,0}u(x)=0\}, \, \qquad {\rm if}\ b=0,\\
 D(\Ab)&=&\{ u\in C^1([0,\infty[)\cap C^2(]0,\infty[)\cap C([0,\infty])\,\mid\, \\
& &\lim_{x\to 0^+} xu''(x)=0,\ \lim_{x\to + \infty} \Ab u(x)=0 \}, \qquad {\rm if}\ b>0.
\end{eqnarray*}
It is known that $(\Ab, D(\Ab))$ generates a $C_0$-contractive  semigroup in $C([0,\infty])$ (see, e.g., \cite{CT,T}). 

\begin{prop}\label{p:generator} Let $b\geq 0$, $\gamma>0$ and let $D$ be  defined according to \eqref{eq:core}. Then the following properties hold.
\begin{itemize} 
\item[\rm (1)] $D(\Ab)\subseteq \{ u\in C^2(]0,\infty[)\,\mid\, \lim_{x\to 0}x u'(x)=0, \ \lim_{x\to\infty}u'(x)=0\}$.
\item[\rm (2)] $D$ is a core for $(A^{\gamma,b}, D(A^{\gamma,b}))$.
\end{itemize}
\end{prop}

{\sc Proof.}  (1) Let $u\in D(\Ab)$. Then, 
for every $\varepsilon >0$ there exists $M>0$ such that  $|A^{\gamma,b}u(x)|<\gamma\varepsilon$ for every $x>M$. On the other hand, we have, for every $x>M$, that 
\begin{eqnarray*}  
\int_M^x A^{\gamma, b} u(s)ds &=& \gamma x u'(x)-\gamma Mu'(M) + \int_M^x (b-\gamma) u'(s)ds  \\
&=&\gamma x u'(x)-\gamma Mu'(M) + (b-\gamma) (u(x)-u(M)) 
\end{eqnarray*}
and hence,
\[
u'(x)=\frac{1}{\gamma x}\int_M^x A^{\gamma, b} u(s)ds+\frac{M}{x}u'(M)+\frac{b-\gamma}{\gamma x}(u(M)-u(x)).
\]
So, we deduce, for every $x>M$, that
\begin{eqnarray*}
 |u'(x)| &\leq& \frac{1}{\gamma x} \int_M^x |A^{\gamma, b} u(s)| ds + \frac{1}{\gamma x}\left[\gamma M|u'(M)| + |(b-\gamma) (u(M)-u(x)|\right]\\ 
&\leq& \varepsilon  +  \frac{1}{\gamma x}\left[\gamma M|u'(M)| + |(b-\gamma) (u(M)-u(x))|\right]. 
\end{eqnarray*}
This ensures  that $\limsup_{x\to+\infty} |u'(x)|\leq \varepsilon$; it follows that $\lim_{x\to+\infty}u'(x)=0$ as $\varepsilon$ is arbitrary.

If $b>0$, then clearly $\lim_{x\to 0^+}xu'(x)=0$. If $b=0$,   we can apply the same argument as before integrating on a suitable small interval $[0,\delta]$.

(2) Let $u\in D(\Ab)$ and let $\phi\in C_c([0,\infty[)$ satisfy $\phi(x)=0$ if $x>2$ and $\phi(x)=1$ if $x<1$.
 For each  $n\in\N$ we set
\[ 
u_n(x):=\left\{ \begin{array}{ll}u(\frac 1 n) + (x-\frac 1 n) u'(\frac 1 n)+ \frac 1 2 (x-\frac 1 n)^2 u''(\frac 1 n) \qquad &0 \leq x \leq \frac 1 n\\
u(x) \qquad &  \frac 1 n\leq x\leq 1\\
(u(x)-l)\phi\left(\frac x n\right)+l \qquad &x\geq 1, 
\end{array} \right.
\]
where $l:=\lim_{x\to \infty} u(x)$.
Then  $u_n\in C^2([0,\infty[)\cap C([0,\infty])$ and $u_n$ is constant in a neighborhood of $+\infty$ as it is easy to verify. So, $(u_n)_n\subset D$.
It is also straightforward to prove that $u_n \rightarrow u$ uniformly on $[1,+\infty[$. 
Moreover, for every $n\in\N$, we have 
\[ \sup_{x\in [0,1]} |u_n(x)-u(x)|\leq \sup_{0\leq x\leq \frac 1 n} \left|u(x)-u\left(\frac 1 n\right)\right| + \frac 1 n \left|u'\left(\frac 1 n\right)\right| +  \frac{1}{n^2} \left|u''\left(\frac 1 n\right)\right|. \]
Since  $\lim_{x\to 0^+}xu'(x)=0$ by (1), it follows that $u_n\rightarrow u$ uniformly on $[0,1]$.

Therefore, $u_n\rightarrow u$ uniformly on $[0,\infty[$.

On the other hand,  for every $n\in\N$, we have 
\begin{eqnarray*}
& & \sup_{x\in [0,1]} |\Ab u_n(x)-\Ab u(x)|= \\
&=& \sup_{0\leq x\leq \frac 1 n} \left|\gamma x u''\left(\frac 1 n\right) + b\left[u'\left(\frac 1 n\right) + \left(x-\frac 1 n\right) u''\left(\frac 1 n\right)\right]-\gamma xu''(x)-bu'(x)\right| \\
&\leq& \sup_{0\leq x\leq \frac 1 n} \left(|\gamma xu''(x)|+ \frac{ \gamma +b}{n} \left| u''\left(\frac 1 n\right)\right| + b\left|u'(x)-u'\left(\frac 1 n\right)\right| \right)
\end{eqnarray*}
(observe that the term $b\left|u'(x)-u'\left(\frac 1 n\right)\right| $ disappears in the above inequality for $b=0$) and 
\begin{eqnarray*}
& &\sup_{x\geq 1}|A^{\gamma, b}u_n(x) -A^{\gamma, b}u(x)|=\\
&= &\sup_{x\geq 1}\left| \left( \phi\left( \frac x n\right) -1\right) A^{\gamma,b}u (x) + 2 \frac{\gamma x}{n} u'(x) \phi'\left( \frac x n\right) + 
\left[\frac{\gamma x}{n^2}\phi''\left( \frac x n\right)+\frac{b}{n}\phi'\left( \frac x n\right)\right](u(x) -l)\right| \\
&\leq& \sup_{x\geq n} \left| \left( \phi\left( \frac x n\right) -1\right) A_{\gamma,b}u (x)\right| \\
& & + \sup_{n\leq x\leq 2n}\left| \frac{2\gamma x}{ n} u'(x) \phi'\left( \frac x n\right) + 
 \left[\frac{\gamma x}{n^2}\phi''\left( \frac x n\right)+\frac{b}{n}\phi'\left( \frac x n\right)\right](u(x)-l) \right| \\
&\leq& \sup_{x\geq n} \left|  A^{\gamma,b}u (x)\right| + 4||\phi'||_\infty\gamma \sup_{n\leq x\leq 2n} |u'(x)| + 
 \left(\frac{2\gamma}{n}||\phi''||_\infty+ \frac{b}{n}||\phi'||_\infty\right)\sup_{n\leq x\leq 2n}|u(x)-l| .
\end{eqnarray*}
Taking (1) into account, it follows that $\Ab u_n\rightarrow \Ab u$ uniformly on $[0,+\infty[$. 

Therefore, $D$ is a core for $(A^{\gamma,b}, D(A^{\gamma,b}))$. \qed

\begin{remark} \rm For similar results in a  more general setting we refer to \cite{ALM}.
\end{remark}

\begin{prop} Let $b\geq 0$ and $\gamma>0$. Then $(\Ab, D(\Ab))$ is the infinitesimal generator of $(\Pb)_{t \geq 0}$. 
\end{prop}

{\sc Proof.} By  Propositions \ref{p:c0semig} and \ref{p:generator}(2) we have that $D\subseteq D(A)$,  $A$ and $A^{\gamma,b}$ coincide on $D$ and that $(\lambda I - A)(D)= (\lambda I - A^{\gamma, b})(D)$ is dense in $C([0,\infty])$ for some $\lambda>0$. It follows that $D$ is also a core for $(A, D(A))$. Since $(A,D(A))$ and $(A^{\gamma,b}, D(A^{\gamma,b}))$ are closed operators in $C([0,\infty])$, the thesis follows. 

\medskip

A straightforward application of the First Trotter-Kato Approximation theorem (see, e.g., \cite[Chap. III, Theorem 4.8]{EN}) and of the results above gives the following result.

\begin{corollary} Let $b\geq 0$ and $\gamma>0$. Then, for every $f\in C([0,\infty])$,
\[ \lim_{b\to 0^+} \Pb f= \Po f\]
in $C([0,\infty])$ uniformly for $t$ in compact intervals.
\end{corollary} 

\section{Analyticity constants for  $ (\Pb)_{t\geq 0}$ in $C([0,\infty])$}

It is known that, for every $b\geq 0$ and $\gamma>0$, the operator $(\Ab, D(\Ab))$
generates a bounded analytic $C_0$--semigroup of angle $\pi/2$ in $C([0,\infty])$. Indeed, $\Ab u$ has the same behaviour in $0$ of the operator
$B^{\gamma,b}$, defined by $ B^{\gamma,b}u=xu''(x) + bu'(x)$, 
with domain 
\begin{eqnarray*}& &D(B^{\gamma,0})=\{ u\in C([0,1]) \cap C^2(]0,1[)\, \mid\, \lim_{x\to 0} Bu =0,\ u'(1)=0\}\qquad if\ b=0,\\
& & D(B^{\gamma,b})=\{ u\in C^1([0,1])\cap C^2(]0,1[)\,\mid\, \lim_{x\to 0} xu''(x)=0,\ u'(1)=0 \}\qquad if\ b>0,
\end{eqnarray*}
which generates a bounded analytic $C_0$--semigroup of angle $\pi/2$ in $C([0,1])$, see \cite{Met,CM}.
On the other hand, by performing the change of variable $x=\frac 1 y$, it can be seen that $\Ab$ behaves at $\infty$ as the operator
$C^{\gamma,b}$, defined by $C^{\gamma,b}v=\gamma y^3v''+(2-b)y^2v'$, with domain 
\[ 
D(C^{\gamma, b})= \{ u\in C([0,1]) \cap C^2(]0,1[)\, \mid\, \lim_{x\to 0^+} C^{\gamma,b}u =0,\ u'(1)=0\}
\]
behaves near $0$. By \cite[Theorem 4.20]{CMPR} and the comments below,  $(C^{\gamma,b}, D(C^{\gamma, b}))$ generates a bounded analytic $C_0$--semigroup of angle $\pi/2$ in  $C([0,1])$.

Now, by suitable cut and paste tecniques (see, e.g., \cite[Proposition 2.4]{CM}), it easily follows that $\Ab$ generates a  bounded analytic $C_0$--semigroup of angle $\pi/2$ in $C([0,\infty])$.
So, there exists $M=M(\gamma,b)>0$ such that $||tA^{\gamma, b}\Pb  ||\leq M_{\gamma, b}$ for every $t\geq 0$.
Nevertheless, we  now show  that $M(\gamma, b)$ is uniformly bounded in bounded intervals  $[0,B]$ and in half--lines $[\gamma_0,\infty[$ with $B,\gamma_0>0$.

\begin{prop} Let $B,\gamma_0>0$.
Then, for every  $b\in [0,B]$,   $\gamma\geq \gamma_0$ and $f\in C([0,\infty])$,  
\begin{equation}\label{eq:tAP} 
||tA^{\gamma,b}\Pb f||_\infty \leq \frac{2(1+\sqrt{2}+b)}{\gamma}||f||_\infty \leq  \frac{2(1+\sqrt{2}+B)}{\gamma_0}||f||_\infty
,\qquad t\geq 0. 
\end{equation}
\end{prop}

{\sc Proof.} Fix   $b>0$ and $\gamma>0$. 
If  $f\in C([0,\infty])$, then straightforward calculations
(see, e.g., \cite[Lemma 4.5]{BP}) show, for every $x\geq 0$ and $t>0$, that  
\begin{eqnarray}\label{eq:derb}
\label{eq:nderb} & & (\Pb f)'(x)= \frac{e^{-\frac{x}{\gamma t}}}{\gamma t} \sum_{m=0}^\infty \left(\frac{x}{\gamma t}\right)^m \frac{1}{m!} \int_0^\infty e^{-z}\left(\frac{z^{m+\frac b \gamma}}{\Gamma (m+\frac b \gamma+1)}-\frac{z^{m+\frac b \gamma -1}}{\Gamma (m+\frac b \gamma)}\right)  \nonumber\\
& & \qquad\qquad\qquad\qquad\times f(\gamma z t)dz,
\end{eqnarray}
\begin{eqnarray}
\label{eq:dersb} &  &(\Pb f)''(x)= \frac{1}{\gamma t} e^{-\frac{x}{\gamma t}}\sum_{m=1}^\infty \left(\frac{x}{\gamma t}\right)^{m-1} \frac{1}{(m-1)!}\\
& &\times \int_0^\infty e^{-z}
\left(\frac{z^{m+\frac b \gamma}}{\Gamma (m+\frac b \gamma+1)}-2\frac{z^{m+\frac b \gamma -1}}{\Gamma (m+\frac b \gamma)} + \frac{z^{m+\frac b \gamma -2}}{\Gamma (m+\frac b \gamma-1)}\right) f(\gamma z t)\frac{dz}{\gamma t}.\nonumber
\end{eqnarray}
So, by Lemma \ref{le:gamma} we obtain, for every $x\geq 0$ and $t>0$, that 
\begin{eqnarray}
& & |(\Pb f)'(x)|\leq ||f||_\infty \frac{1}{\gamma t} e^{-\frac{x}{\gamma t}}\sum_{m=0}^\infty \left(\frac{x}{\gamma t}\right)^m \frac{1}{m!\Gamma(m + \frac b \gamma +1)}\times \nonumber\\
& & \times \int_0^\infty e^{-z}z^{m+\frac b \gamma -1}\left|z-m-\frac b \gamma\right|dz \nonumber\\
&\leq &||f||_\infty \frac{1}{\gamma t} e^{-\frac{x}{\gamma t}}\sum_{m=0}^\infty \left(\frac{x}{\gamma t}\right)^m \frac{2}{m!\sqrt{m+\frac b \gamma}}\leq  ||f||_\infty \frac{2}{\gamma t}. \label{eq:4}
\end{eqnarray}
On the other hand, summing by parts in \eqref{eq:dersb} we have, for every $x\geq 0$ and $t>0$, that  
\begin{eqnarray*}
& & x(\Pb f)''(x) =-\frac{1}{\gamma t}e^{-\frac{x}{\gamma t}} \frac{x}{\gamma t} \int_0^\infty e^{-z} f(\gamma zt) \left[ \frac{z^{b/\gamma}}{\Gamma(\frac b \gamma +1)} - \frac{z^{b/\gamma-1}}{\Gamma(\frac b \gamma)}\right] dz +\\ 
& &+\frac{1}{\gamma t}e^{-\frac{x}{\gamma t}} \sum_{m=1}^\infty \left( \frac{x}{\gamma t} \right)^m \frac{1}{m!} \left(m-\frac{x}{\gamma t}\right) 
\int_0^\infty  e^{-z}f(\gamma zt) \left[ \frac{z^{m+b/\gamma}}{\Gamma(m+\frac b \gamma +1)} - \frac{z^{m+b/\gamma-1}}{\Gamma(m+\frac b \gamma)}\right] dz.
\end{eqnarray*}
Therefore, again by Lemma \ref{le:gamma} we deduce,  for every $x\geq 0$ and $t>0$, that
\begin{eqnarray}\label{eq:4-1}
& & |x(\Pb f)''(x)| \leq \frac{2}{\gamma t}||f||_\infty + \nonumber\\
& & +\frac{1}{\gamma t}||f||_\infty e^{-\frac{x}{\gamma t}} \sum_{m=1}^\infty \left( \frac{x}{\gamma t} \right)^m \frac{1}{m!} \left|m-\frac{x}{\gamma t}\right| 
\int_0^\infty\frac{z^{m+b/\gamma-1}}{\Gamma(m+\frac b \gamma +1)}\left|z-m-\frac b \gamma\right| dz\nonumber\\
& &\leq \frac{2}{\gamma t}||f||_\infty+\frac{1}{\gamma t}||f||_\infty e^{-\frac{x}{\gamma t}} \sum_{m=1}^\infty \left( \frac{x}{\gamma t} \right)^m \frac{1}{m!} \left|m-\frac{x}{\gamma t}\right| \frac{2}{\sqrt{m}}\nonumber\\
& & \leq  \frac{2(1+\sqrt{2})}{\gamma t}||f||_\infty,
\end{eqnarray}
as by H\"older inequality we have 
\begin{eqnarray}\label{eq:holder}
\sum_{m=1}^\infty \frac{s^m}{m!}|s-m| m^{-\frac 1 2} &\leq & \left( \sum_{m=1}^\infty \frac{s^m}{m!}[s-m]^2\right)^\frac 1 2 \left( \sum_{m=1}^\infty \frac{s^m}{m!m}\right)^\frac 1 2\nonumber\\
& \leq &\sqrt{e^s s}\sqrt{\frac{2e^s}{s}}= \sqrt{2} e^s,\quad s>0.
\end{eqnarray}

Now, let $b=0$. If  $f\in C([0,\infty])$, then straightforward calculations
(see, e.g., \cite[Lemma 4.1]{BP}) show, for every $x\geq 0$ and $t>0$, that  
\begin{eqnarray}
\label{eq:nder}(P_t^{\gamma,0}f)'(x)&=& e^{-\frac{x}{\gamma t}}\int_0^\infty [f(\gamma zt)-f(0)]e^{-z}\frac{dz}{\gamma t}\\
& & +e^{-\frac{x}{\gamma t}}\sum_{m=1}^\infty\left(\frac{x}{\gamma t}\right)^m\frac{1}{m!}\int_0^\infty f(\gamma zt)e^{-z}\left[\frac{z^m}{m!}-\frac{z^{m-1}}{(m-1)!}\right]\frac{dz}{\gamma t},\nonumber\\
\label{eq:nders}(P_t^{\gamma,0} f)''(x)&=&\frac{1}{\gamma t}e^{-\frac{x}{\gamma t}}\int_0^\infty [f(\gamma zt)-f(0)]e^{-z}(z-2)\frac{dz}{\gamma t}\\
& &+\frac{1}{\gamma t}e^{-\frac{x}{\gamma t}}\sum_{m=1}^\infty\left(\frac{x}{\gamma t}\right)^m\frac{1}{m!}\int_0^\infty f(\gamma zt)e^{-z}\left[\frac{z^{m+1}}{(m+1)!}-2\frac{z^m}{m!}+\frac{z^{m-1}}{(m-1)!}\right]\frac{dz}{\gamma t}.\nonumber
\end{eqnarray}

Since   $P_t^{\gamma,0}g=P_t^{\gamma,0} f-f(0)$ if $g=f-f(0)$ and so $(P_t^{\gamma,0}g)'= (P_t^{\gamma,0} f)'$ and $(P_t^{\gamma,0}g)''= (P_t^{\gamma,0} f)''$, w.l.o.g. we may suppose $f(0)=0$. Therefore, summing by parts in \eqref{eq:nders} we  have, for every $x\geq 0$ and $t>0$, that
\begin{eqnarray*}
& & x(P_t^{\gamma,0} f)''(x)= -e^{-\frac{x}{\gamma t}}\frac{x}{\gamma t}\frac{1}{\gamma t}\int_0^\infty f(\gamma z t)e^{-z}dz+\\
&  & + e^{-\frac{x}{\gamma t}}\frac{1}{\gamma t}\sum_{m=1}^\infty \left(\frac{x}{\gamma t}\right)^m\frac{1}{m!}\left(\frac{x}{\gamma t}-m\right)\int_0^\infty f(\gamma z t)e^{-z}\frac{z^{m-1}}{m!}(m-z)dz.
\end{eqnarray*}
By Lemma \ref{le:gamma} and \eqref{eq:holder} we obtain,  for every $x\geq 0$ and $t>0$, that 
\begin{eqnarray}\label{eq:43}
& & | x(P_t^{\gamma,0} f)''(x)|\leq \frac{||f||_\infty}{\gamma t}e^{-\frac{x}{\gamma t}}\frac{x}{\gamma t}\nonumber\\
& &+ e^{-\frac{x}{\gamma t}}\frac{1}{\gamma t}\sum_{m=1}^\infty \left(\frac{x}{\gamma t}\right)^m\frac{1}{m!}\left|\frac{x}{\gamma t}-m\right|\int_0^\infty f(\gamma z t)e^{-z}\frac{z^{m-1}}{m!}|m-z|dz\nonumber\\
& & \leq \frac{||f||_\infty}{\gamma t}e^{-\frac{x}{\gamma t}}\left[\frac{x}{\gamma t}+2\sum_{m=1}^\infty \left(\frac{x}{\gamma t}\right)^m\frac{1}{m!\sqrt{m}}\left|\frac{x}{\gamma t}-m\right|\right]\nonumber\\
& & \leq \frac{||f||_\infty}{\gamma t}e^{-\frac{x}{\gamma t}}\left[\frac{x}{\gamma t}+2\sqrt{2}e^{-\frac{x}{\gamma t}}\right]\leq \frac{1+2\sqrt{2}}{\gamma t}||f||_\infty.
\end{eqnarray}

Combining \eqref{eq:4} and \eqref{eq:4-1} (by \eqref{eq:43} in case $b=0$), it follows that
\[
 ||tA^{\gamma, b}\Pb f||_\infty \leq \frac{2(1+\sqrt{2}+b)}{\gamma} ||f||_\infty\leq  \frac{2(1+\sqrt{2}+B)}{\gamma_0}||f||_\infty, 
\]
for every $ t\geq 0,\, b\in [0,B], \gamma\in [\gamma_0,\infty[$ and $f\in C([0,\infty])$.\qed

\medskip

As a consequence, we obtain  (see, e.g.,    \cite[Chap. II, Theorem 4.6]{EN}). 

\begin{corollary}\label{coro:res} Let $B,\gamma_0>0$. Then
there exists $\theta_1=\theta(B, \gamma_0) \in ]0, \frac\pi 2]$ such that, 
 for every $\theta \in ]0,\theta_1[$ there exists  $d_1=d(B, \gamma_0,\theta )>0$  for which
\begin{equation}\label{eq:47} 
||R(\lambda,\Ab) f||_\infty \leq d_1 \frac{||f||_\infty}{|\lambda|},\ \ |{\rm arg}\lambda|\leq \frac \pi 2 +\theta,\, \lambda\not=0,\end{equation}
for every $f\in C([0,\infty])$, $b\in [0,B]$ and $\gamma\geq \gamma_0$.
\end{corollary}

\section{Gradient estimates}

\begin{prop}\label{prop:derres} Let $ B,\gamma_0>0$. Then, for every $b\in [0,B]$ and  $\gamma\geq \gamma_0$,     
 the following properties hold.
 \begin{itemize} 
 \item[\rm (1)] For every  $f\in C([0,\infty])$ and  $t>0$, $\Pb f\in C^1([0,\infty[)$ and 
 \begin{equation}\label{eq:5}
 |\sqrt{x}(\Pb f)'(x)| \leq \frac{2\sqrt{2}}{\sqrt{\gamma_0 t}}||f||_\infty,\quad x\geq 0.
\end{equation} 
 \item[\rm (2)] There exists $\theta_1=\theta(B, \gamma_0) \in ]0, \frac\pi 2]$ such  that, 
 for every $\theta \in ]0,\theta_1[$ there exists  $d_2=d(B, \gamma_0,\theta)>0$  for which 
\begin{eqnarray}
& & \label{eq:limres0-1} R(\lambda,\Ab) f \in C^1(]0,\infty[),\\
& &\label{eq:limres0} \lim_{x\to 0^+} \sqrt{x}(R(\lambda, \Ab)f)'(x)=0,\\
& & \label{eq:ressec0}|\sqrt{x}(R(\lambda,\Ab) f)'(x)| \leq d_2 \frac{||f||_\infty}{\sqrt{|\lambda|}},  \quad x>0, 
\end{eqnarray}
for every $f\in C([0,\infty])$ and $|{\rm arg}\lambda|\leq \frac \pi 2 +\theta$ with $\lambda\not=0$, where $d_1$ is the constant appearing in \eqref{eq:47}.
\end{itemize}
\end{prop}
 
 {\sc Proof.}  
(1) Let $b>0$ and $\gamma>0$. From  \eqref{eq:nderb} it follows immediately that $\Pb f\in C^1([0,\infty[)$. 
 On the other hand, we have
\begin{eqnarray*}
\frac{1}{\Gamma(\frac{b}{\gamma}+1)}\int_0^\infty e^{-z}z^{\frac{b}{\gamma}-1}\left|z-\frac{b}{\gamma}\right|dz
&  \leq &\frac{1}{\Gamma(\frac{b}{\gamma}+1)}\int_0^\infty e^{-z}z^{\frac{b}{\gamma}}dz+\frac{\frac{b}{\gamma}}{\Gamma(\frac{b}{\gamma}+1)}\int_0^\infty e^{-z}z^{\frac{b}{\gamma}-1}dz\\
&  = & 1+\frac{\frac{b}{\gamma}}{\Gamma(\frac{b}{\gamma}+1)}\Gamma(\frac{b}{\gamma})=2.
\end{eqnarray*}
So, fixed any $f\in C([0,\infty])$ and 
 applying H\"older inequality in \eqref{eq:4} and Lemma \ref{le:gamma}, we get, for every $t>0$ and $x>0$, that
\begin{eqnarray*}
 |(\Pb f)'(x)|&\leq& ||f||_\infty \frac{1}{\gamma t} e^{-\frac{x}{\gamma t}}\sum_{m=0}^\infty \left(\frac{x}{\gamma t}\right)^m \frac{1}{m!\Gamma(m + \frac b \gamma +1)} \int_0^\infty e^{-z}z^{m+\frac b \gamma -1}\left|z-m-\frac b \gamma\right|dz\\
&\leq&||f||_\infty \frac{1}{\gamma t} e^{-\frac{x}{\gamma t}}\left(2+\sum_{m=1}^\infty \left(\frac{x}{\gamma t}\right)^m \frac{1}{m!\sqrt{m+\frac b \gamma}}\right)\\
&\leq &||f||_\infty \frac{1}{\gamma t} e^{-\frac{x}{\gamma t}}\left(2+2\sum_{m=1}^\infty \left(\frac{x}{\gamma t}\right)^m \frac{1}{m!\sqrt{m+1}}\right)\\
&=& 2||f||_\infty \frac{1}{\gamma t} e^{-\frac{x}{\gamma t}}\sum_{m=0}^\infty\left(\frac{x}{\gamma t}\right)^m \frac{1}{m!\sqrt{m+1}}\\
 &\leq& 2||f||_\infty \frac{1}{\gamma t} e^{-\frac{x}{\gamma t}}\left(\sum_{m=0}^\infty \left(\frac{x}{\gamma t}\right)^m \frac{1}{m!(m+1)}\right)^{\frac 1 2} e^{\frac{x}{2\gamma t}}\\
 &=& 2||f||_\infty \frac{1}{\gamma t} e^{-\frac{x}{2\gamma t}}\left[\frac{\gamma t}{x} (e^{\frac{x}{\gamma t}}-1)\right]^{\frac 1 2}\\
 &=&  2||f||_\infty \frac{1}{\sqrt{\gamma t x}}\left( 1- e^{-\frac{x}{2\gamma t}}\right)^\frac 1 2 \leq 2||f||_\infty \frac{1}{\sqrt{\gamma t x} }.
 \end{eqnarray*}
Now, let $b=0$. Then by \eqref{eq:nder} we have, for every $t>0$ and $x>0$, that
\begin{eqnarray*}
& & |(P_t^{\gamma,0} f)'(x)|\leq 2 ||f||_\infty \frac{e^{-\frac{x}{\gamma t}}}{\gamma t} \left(1+\sum_{m=1}^\infty\left(\frac{x}{\gamma t}\right)^m \frac{1}{m!\sqrt{m}}\right)\\
& &\leq  2||f||_\infty \frac{e^{-\frac{x}{\gamma t}}}{\gamma t} \left(\sqrt{2}+\sum_{m=1}^\infty\left(\frac{x}{\gamma t}\right)^m \frac{\sqrt{2}}{m!\sqrt{m+1}}\right)\\
& &= 2\sqrt{2}||f||_\infty \frac{e^{-\frac{x}{\gamma t}}}{\gamma t}\sum_{m=0}^\infty\left(\frac{x}{\gamma t}\right)^m \frac{1}{m!\sqrt{m+1}}
\leq  2\sqrt{2}||f||_\infty\frac{1}{\sqrt{\gamma t x} }.
\end{eqnarray*}
So, (1) is satisfied.

(2) By (1) and the dominated convergence theorem we obtain, for every $\eta>1$ and $x>0$, that
\[
\sqrt{x}(R(\eta,\Ab)f)'(x)=\sqrt{x}D \left(\int_0^{+\infty} e^{-\eta t} \Pb fdt\right)=\int_0^\infty e^{-\eta t} \sqrt{x} (\Pb f)'dt.\]
 At this point,  \eqref{eq:limres0} and \eqref{eq:ressec0} follow  by arguing as in the proof of \cite[Proposition 2.1]{AM}.\qed

\medskip

In case   $b>0$,  we can achieve the following global gradient estimate for the resolvent.

\begin{prop}\label{res} Let $b,\gamma>0$. Then 
there exists a constant $C>0$ independent on $b$ such that, for every $\lambda>0$ and $f\in C([0,\infty])$, we have
\[
 ||(R(\lambda, A_b)f)'||_\infty \leq \frac{1}{\gamma}\max\left\{2,\frac{\gamma}{b}\right\} C ||f||_\infty .
\]
\end{prop} 

In order to prove this, we need the following result.

\begin{lemma} \label{l.derivata}
Let $f\in C^1([0,\infty[) \cap C([0,\infty])$, with $f'$ bounded. Then, for  every $x\geq 0$, we have
\[
 (\Pb f)'(x) = e^{-\frac{ x}{\gamma t}} \sum_{m=0}^\infty \left(\frac{x}{\gamma t}\right)^m \frac{1}{m!\Gamma(m+\frac b \gamma +1)} \int_0^\infty f'(\gamma zt)z^{m+\frac b \gamma}e^{-z}dz.
\] 
\end{lemma}

{\sc Proof.} 
Integrating by parts in \eqref{eq:nderb}, we obtain, for every $x\geq 0$ and $t>0$, that
\begin{eqnarray*}
& &(\Pb f)'(x) = \frac{e^{-\frac{x}{\gamma t} }}{\gamma t} \sum_{m=0}^\infty \left( \frac {x}{\gamma t}\right)^m \frac{1}{m!} \int_0^\infty e^{-z}f(\gamma zt)\left( \frac{z^{m+b/\gamma }}{\Gamma(m+b/\gamma +1)} - \frac{z^{m+b/\gamma-1}}{\Gamma(m+b/\gamma)}\right) dz\\
& & =\frac{e^{-\frac{x}{\gamma t} }}{\gamma t} \sum_{m=0}^\infty \left( \frac{x}{\gamma t}\right)^m \frac{1}{m!}\left( \int_0^\infty e^{-z}f(\gamma zt) \frac{z^{m+b/\gamma}}{\Gamma(m+b/\gamma+1)}dz - \int_0^\infty e^{-z}f(\gamma z t) \frac{z^{m+b/\gamma-1}}{\Gamma(m+b/\gamma)}dz\right)\\
& & = \frac{e^{-\frac{x}{\gamma t} }}{\gamma t} \sum_{m=0}^\infty \left( \frac{x}{\gamma t}\right)^m \frac{1}{m!} 
\left( \left[ -e^{-z}f(\gamma zt) \frac{z^{m+b/\gamma}}{\Gamma(m+b/\gamma+1)} \right]_0^\infty  + \int_0^\infty\gamma t e^{-z} f'(\gamma zt)\frac{z^{m+b/\gamma}}{\Gamma(m+b/\gamma +1)} dz \right. \\
& &\ \left. +\int_0^\infty e^{-z}f(\gamma zt) (m+b/\gamma) \frac{z^{m+b/\gamma-1}}{\Gamma(m+b/\gamma+1)} dz - 
\int_0^\infty e^{-z}f(\gamma zt)  \frac{z^{m+b/\gamma-1}}{\Gamma(m+b/\gamma)}dz\right)\\
&=& e^{-\frac{ x}{\gamma t}} \sum_{m=0}^\infty \left(\frac{x}{\gamma t}\right)^m \frac{1}{m! \Gamma(m+b/\gamma +1) } \int_0^\infty f'(\gamma zt)z^{m+b/\gamma }e^{-z}dz. \qquad \qed
\end{eqnarray*}

We observe that the above inequality ensures that, for every  $f\in C^1([0,\infty[) \cap C([0,\infty])$, with $f'$ bounded,   we have
\[ 
||(R(\lambda, \Ab)f)'||_\infty \leq ||f'||_\infty,\quad \lambda >0.\qquad\
\]
But in order  to estimate $||(R(\lambda, \Ab)f)'||_\infty$ with $||f||_\infty$, we need to proceed as follows.

\medskip

{\sc Proof of Proposition \ref{res}.} We first assume that 
$f\in C^1([0,\infty[)\cap C([0,\infty])$,  with $f'$ bounded. In such a case, by  Lemma \ref{l.derivata} we have, for every $x\geq 0$ and $\lambda>0$, that 
\begin{eqnarray*}
& &(R(\lambda, \Ab)f)'(x)= \int_0^\infty e^{-\lambda t}(\Pb f)'(x)dt \\
& & =\int_0^\infty e^{-\lambda t} e^{-\frac{x}{\gamma t}}\left( \sum_{m=0}^\infty \left( \frac{x}{\gamma t}\right)^m \frac{1}{m!} \frac{1}{\Gamma(m+\frac b \gamma +1)} \int_0^\infty f'(\gamma z t) z^{m+\frac b \gamma} e^{-z}dz\right) dt \\
& & = \sum_{m=0}^\infty \int_0^\infty dz \left(e^{-z} z^{m+\frac b\gamma} \frac{1}{m!} \frac{1}{\Gamma(m+\frac b\gamma +1)} \int_0^\infty e^{-\lambda t} \left( \frac{x}{\gamma t}\right)^m f'(\gamma z t)e^{-\frac{x}{\gamma t}}dt\right).
\end{eqnarray*}
 On the other hand, for every $m\geq 0$, $x\geq 0$ and $\lambda>0$, we have
\begin{eqnarray*}
& &\int_0^\infty e^{-\lambda t} \left( \frac{x}{\gamma t}\right)^m f'(\gamma z t)e^{-\frac{x}{\gamma t}}dt=\\
& & =\left[ \frac{f(\gamma z t)}{ \gamma z} e^{-\lambda t} \left(\frac{x}{\gamma t}\right)^m e^{-\frac{x}{\gamma t}} \right]_0^\infty 
- \frac{1}{\gamma z} \int_0^\infty f(\gamma zt) \left[ e^{-\lambda t} \left(\frac{x}{\gamma t}\right)^m e^{-\frac{x}{\gamma t}}\right]'dt\\
& & =\frac{1}{\gamma z} \int_0^\infty e^{-\lambda t} \left(\frac{x}{\gamma t}\right)^m e^{-\frac{x}{\gamma t}}\left(\lambda + \frac{m}{t} - \frac{x}{\gamma t^2}\right) f(z\gamma t)dt.
\end{eqnarray*}
So, we obtain, for every $x\geq 0$ and $\lambda>0$, that 
\begin{eqnarray*} 
& &(R(\lambda, \Ab)f)'(x)=\\
& & =\frac 1 \gamma \sum_{m=0}^\infty \frac{1}{m!\Gamma (m+\frac b\gamma +1)}\int_0^\infty dz e^{-z} z^{m+\frac b \gamma -1} \int_0^\infty f(\gamma zt) e^{-\lambda t}\left( \frac{x}{\gamma t}\right)^m \times\\
&  &\times e^{-\frac{x}{\gamma t}}\left(\lambda + \frac{m}{t} - \frac{x}{\gamma t^2}\right) dt \\
& & = \frac 1 \gamma \int_0^\infty dt e^{-\lambda t}\sum_{m=0}^\infty \left( \frac{x}{\gamma t}\right)^m \frac{1}{m!\Gamma (m+\frac b\gamma +1)}e^{-\frac{x}{\gamma t}} 
\left(\lambda + \frac{m}{ t} - \frac{x}{\gamma t^2}\right) \times\\
&  &\times\int_0^\infty e^{-z} z^{m+\frac b \gamma -1}  f(\gamma zt) dz, 
\end{eqnarray*}
 after having observed that we can  interchange sums and integrals because   the  integrals and series are absolutely summable.

 By Lemma \ref{lem:2} it follows, for every $x\geq 0$ and $\lambda>0$, that 
\begin{eqnarray*}
& & |R(\lambda, \Ab)f)'(x)| \leq \frac{||f||_\infty}{\gamma} \int_0^\infty dt e^{-\lambda t}\sum_{m=0}^\infty \left( \frac{x}{\gamma t}\right)^m \frac{1}{m!(m+\frac b \gamma)}e^{-\frac{x}{\gamma t}} 
\left|\lambda + \frac{m}{ t} - \frac{x}{\gamma t^2}\right|\\
& \leq & 
\frac{||f||_\infty}{\gamma} \left( \frac \gamma b \int_0^\infty \lambda e^{-\lambda t}dt + \max\{2,\frac{\gamma}{b}\} \int_0^\infty \frac{e^{-\lambda t}}{ t}
\sum_{m=0}^\infty \left( \frac{x}{\gamma t}\right)^m \frac{1}{(m+1)!}e^{-\frac{x}{\gamma t}} 
\left|m - \frac{x}{\gamma t}\right| dt\right) \\
&\leq & \frac{1}{\gamma}\max\left\{2,\frac{\gamma}{b}\right\}||f||_\infty \left( 1 + \int_0^\infty e^{-\lambda \frac{x}{\gamma s}}\frac{1}{s} \sum_{m=0}^\infty s^m e^{-s}|m-s| \frac{1}{(m+1)!}ds \right) \\
&\leq&  \frac{1}{\gamma}\max\left\{2,\frac{\gamma}{b}\right\} C||f||_\infty 
\end{eqnarray*}
with $C$ independent on $b$. 

Finally,  let $f\in C([0,\infty])$ and let $(f_n)_n \subseteq C^1([0,\infty])$ be an approximating sequence for $f$ in $C([0,\infty])$. Then 
from \eqref{eq:ressec0} it follows that 
\[\lim_{n\to\infty} \sqrt{x}\partial_x R(\lambda,\Ab)f_n(x)=\sqrt{x}\partial_xR(\lambda, \Ab)f(x)
\]
 for every $x\geq 0$ and $\lambda>0$. So, 
 for every $x>0$ and $\lambda>0$,  we obtain
\begin{eqnarray*}
 &  &|(R(\lambda,\Ab)f)'(x)|= \lim_{n\to\infty}|(R(\lambda,\Ab)f_n)' (x)|\\
& & \leq \lim_{n\to\infty} \frac{1}{\gamma}\max\left\{2,\frac{\gamma}{b}\right\}||f_n||_\infty = \frac{1}{\gamma}\max\left\{2,\frac{\gamma}{b}\right\}||f||_\infty. \qquad \qed
\end{eqnarray*}

\begin{corollary}Let $b,\gamma>0$. Then 
there exists a constant $C>0$ independent on $b$ such that, for every $\lambda>0$ and $f\in C([0,\infty])$, we have
\[
 ||x(R(\lambda, A_b)f)''||_\infty \leq  \frac{C}{\gamma}\left[1+\max\left\{2\frac{b}{\gamma},1\right\}\right] ||f||_\infty .
\]
\end{corollary}

{\sc Proof.} Combining  Corollary \ref{coro:res} and Proposition \ref{res} we obtain, for every $\lambda>0$ and $f\in C([0,\infty])$, that
\begin{eqnarray*}
|x(R(\lambda, \Ab)f)''(x)|&=&\left| \frac{1}{\gamma}\left[ \Ab  (R(\lambda,\Ab f) - b(R(\lambda,\Ab)f)'(x)\right]\right|\\   
&= & \left|\frac{1}{\gamma}\left[ \lambda R(\lambda,\Ab)f - f- b(R(\lambda,\Ab)f)'(x)\right]\right| \\
& \leq &\frac{1}{\gamma}\left[1+d_1 +  \frac{b}{\gamma}\max\left\{2,\frac{\gamma}{b}\right\}C \right] ||f||_\infty\\
&\leq & \frac{C'}{\gamma}\left[1+\max\left\{2\frac{b}{\gamma},1\right\}\right] ||f||_\infty,
\end{eqnarray*} 
where $C':=\max\{1+d_1,C\}$. 
 \qed

\section*{Acknowledgements}
The authors thank the anonymous referee for his careful reading of the paper and for many valuable comments and corrections. 

\bibliographystyle{plain}

\end{document}